\documentclass[11pt,a4]{article}

\pdfoutput=1

\usepackage{amsfonts,amsmath,amssymb,bm,graphicx,natbib,hyperref,a4wide}

\begin{document}

\title{\textbf{Numerical methods for the construction of linear waves structures on discrete orthogonal grids}}

\author{Sergey Dvornikov and Maxim Dvornikov\thanks{maxdvo@izmiran.ru} 
\\
\small{\ Pushkov Institute of Terrestrial Magnetism, Ionosphere} \\
\small{and Radiowave Propagation (IZMIRAN),} \\
\small{108840 Troitsk, Moscow, Russia}}

\date{}

\maketitle

\begin{abstract}
We study the possibility for the implementation of linear wave structures on discrete grids with various dimensions. The systems of the first order differential equations for the set of virtual functions, describing the wave propagation, are derived for different spectral zones. We are particularly interested in the spectral regions which are localized near the lowest and the maximal frequencies, as well as in the vicinity of a central frequency. The cases of discrete orthogonal grids with dimensions from 1 to 4 are analyzed. The possibility of the propagation of waves in the considered situations is substantiated by numerical simulations.
\end{abstract}

\section{Introduction}

Studies of relativistic wave equations on discrete grids are of great importance, mainly, for the lattice simulations of quantum fields~\citep{Rat18}. Lattices used in this kind of simulations can have from $\sim 10$ to $\sim 100$ nodes in each dimension. With such a great number of nodes, there are some problems in this research related to the precise approximation of derivatives which one encounters in relativistic wave equations. The relativistic wave equations, such as Dirac, Maxwell, and Klein-Gordon equations, which one typically deals with in quantum field theory~\citep{Wei96}, mainly have the derivatives up to the second order. Therefore, while studying quantum fields on a lattice, one has to derive an expression for the first and the second derivatives of a function defined on a discrete multi-node grid. We should mention that the Dirac equation involves a multi-component wave function. Thus the first derivatives in the corresponding system should be defined consistently.

The problem of numerical differentiation on a multi-node grid was studied by~\cite{Dvo07,Dvo08}, where the expressions for the first and the second derivatives of a function given on a discrete equidistant grid were derived. The spectral properties of the corresponding differentiating filters were analyzed by~\cite{Dvo08}. The results of~\cite{Dvo07,Dvo08} were generalized by~\cite{Mut12,Mut13} to include non-equidistant grids.

In this work, basing on the results of~\cite{Dvo07,Dvo08}, we construct linear wave equations on discrete orthogonal grids of various dimensions. Our work is the generalization of the findings of~\cite{DvoDvo18}, where some of the results were briefly outlined. One of the main goals of the present research is to study the possibility of the existence of the solutions of such wave equations in different spectral regions.

Our work is organized in the following way. We start in Sec.~\ref{sec:DEF} with some definitions. Then, in Secs.~\ref{sec:1D0N/2} and~\ref{sec:1DN/4}, we study linear wave structures in the one dimensional grid in different spectral regions. The two dimensional grid was studied in Secs.~\ref{sec:2D0N/2} and~\ref{sec:2DN/4}. Linear structures in three and four dimensions are analyzed in Secs.~\ref{sec:3D} and~\ref{sec:4D}. We briefly discuss our results in Sec.~\ref{sec:CONCL}.

\section{Definitions}\label{sec:DEF}

Let us recall that, in the construction of linear one dimensional wave structures, at each point $x$ of a one dimensional continuous space, it is necessary to define two differentiable functions $p$ and $q$ satisfying two partial differential equations,
\begin{equation}\label{eq:gensys}
  \frac{\partial p}{\partial \tau} = \frac{\partial q}{\partial x},
  \quad
  \frac{\partial q}{\partial \tau} = \frac{\partial p}{\partial x}.
\end{equation}
The solutions of these equations are harmonic functions of an arbitrary frequency $f_x$  that propagate with the velocity $v=1$  in the positive and negative directions of the coordinate $x$~\citep{TikSam63}: $p=A_1\exp[2\pi \mathrm{i}(f_\tau \tau + f_x x)]$ and $q=A_2\exp[2\pi \mathrm{i}(f_\tau \tau + f_x x)]$, where $f_\tau = \pm f_x$. When constructing wave structures in a two dimensional space, four functions are necessary. In the three dimensional and four dimensional spaces, one needs eight and 16 functions respectively.

In the following, we shall use the one dimensional discrete Fourier transform of the function $S(x)$, given on a discrete grid, with $x= 0,1,\dots,N-1$. It is defined as
\begin{equation*}
  F(f_x) = \sum_{x=0}^{N-1} S(x)
  \exp
  \left[
    -\mathrm{i} \frac{2\pi}{N} f_x\cdot x
  \right],
\end{equation*}
where $-N/2 < f_x \leq N/2$. The two dimensional discrete Fourier transform of the function $S(x,y)$ given on a discrete grid, with $x,y=0,1,\dots,N-1$, has the form,
\begin{equation*}
  F(f_x,f_y) = \sum_{x=0}^{N-1}\sum_{y=0}^{N-1} S(x,y)
  \exp
  \left[
    -\mathrm{i} \frac{2\pi}{N} (f_x\cdot x+f_y\cdot y)
  \right],
\end{equation*}
where $-N/2 < f_{x,y} \leq N/2$. Similar discrete Fourier transforms will be used for three and four dimensional grids.

\section{One dimensional grid: Frequency regions $[0]+[N/2]$}\label{sec:1D0N/2}

\begin{figure}
  \centering
  \includegraphics[scale=0.6]{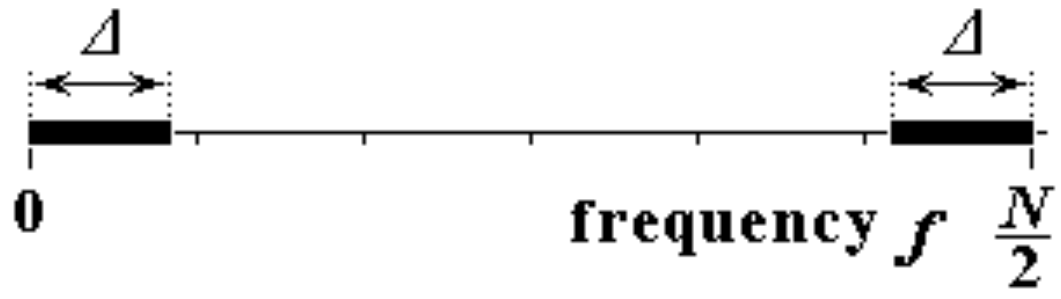}
  \caption{Frequency regions of the discrete function $S(x)$.  
  \label{fig:sch0N/2}}
\end{figure}

Now we consider the one dimensional discrete function $S(x)$, which is defined only for integer values of the argument $x= 0,1,\dots,N-1$. We suppose that the spectrum of this function is located near the zero and maximum spatial frequencies $N/2$, within two frequency intervals: $[0,\Delta]$ and $[N/2–\Delta,N/2]$, as shown in bold segments in Fig.~\ref{fig:sch0N/2}. We shall denote these bands as $[0]+[N/2]$. Using the inverse discrete Fourier transform, one can show that, in this case, $S(x)$ can be represented in the form,
\begin{equation}\label{eq:cos2sin2}
  S(x) = p(x) \cos^2 \frac{\pi x}{2} + q(x) \sin^2 \frac{\pi x}{2},
\end{equation}
where $p$ and $q$ are two linearly independent functions. Hence, for even $x$, $S(x) = p(x)$, whereas for odd $x$ one has $S(x) = q(x)$. Let us call $p$ and $q$ as virtual functions.

\begin{figure}
  \centering
  \includegraphics[scale=0.45]{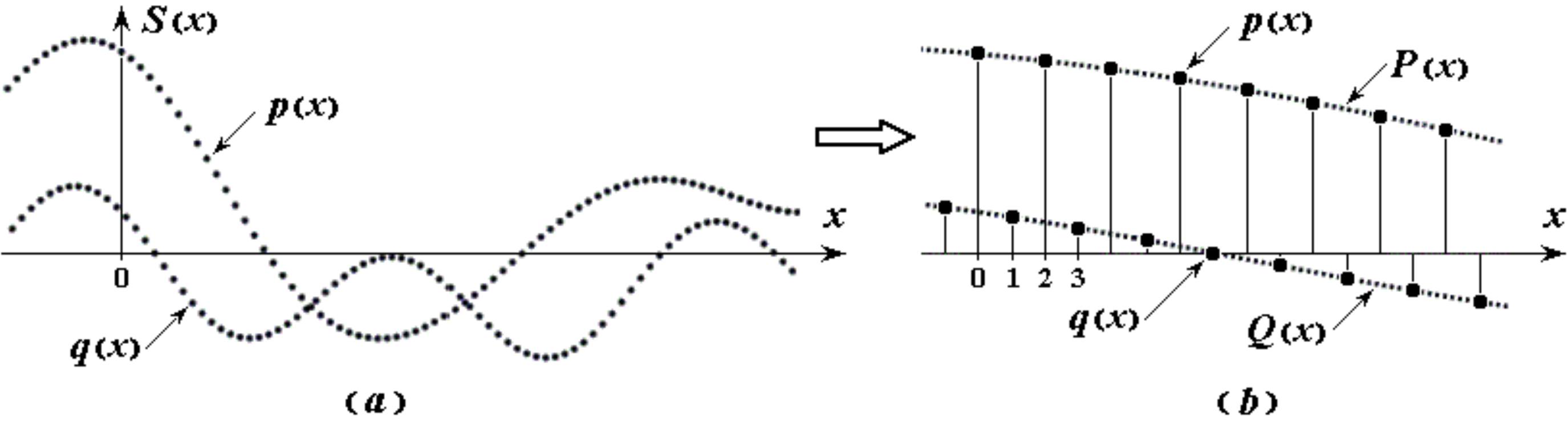}
  \caption{(a) A sample discrete function for the frequency zones $[0]+[N/2]$.
  (b) The fragment of a function.  
  \label{fig:sampfun0N/2}}
\end{figure}

An example of the discrete function $S$ in these frequency bands, provided that $\Delta \ll N$, is given in Fig.~\ref{fig:sampfun0N/2}(a). The fragment of such a function is shown in Fig.~\ref{fig:sampfun0N/2}(b). One can see in these figures that the original discontinuous function $S$ can be considered as the sum of two relatively ``smooth'' functions $p$ and $q$, where $p$ is defined in even $x$, and $q$ in odd $x$.

The differentiation of such ``smooth'' functions is, of course, meaningless by the definition of a derivative. However, if one takes, for example, $2n$ neighboring even values of $S$, i.e. $p(x)$, then one can pass a continuous interpolation polynomial $P(x)$ of the degree $2n - 1$ through these points. In this case, we can take $\mathrm{d}p/\mathrm{d}x \approx \mathrm{d}P/\mathrm{d}x$, i.e. approximate the derivative of the discrete function $p$ with the derivative of the polynomial. The same procedure can be made for the values of $S$ in odd points, i.e. for $q(x)$. In this case one has $\mathrm{d}q/\mathrm{d}x \approx \mathrm{d}Q/\mathrm{d}x$.

Note that the expression for the weight coefficients $\alpha_n(m)$, that allow one to find the first derivative of the discrete function given at even or odd points with any given accuracy, was obtained by~\cite{Dvo07,Dvo08},
\begin{align}\label{eq:firstder}
  f'(x) \approx & \sum_{m=1}^n
  \alpha_n(m)
  [f(x+2m-1)-f(x+2m+1)],
  \notag
  \\
  \alpha_n(m) = &
  \left\{
    2(2m-1)
    \prod_{\substack{k=1 \\ k\neq m}}^{n}
    \left[
      1 - \frac{(2m-1)^2}{(2m-1)^2}
    \right]
  \right\}^{-1},
\end{align}
where $m = 1, 2, \dots, n$. Eq.~\eqref{eq:firstder} was derived by finding the coefficients in the linear terms in $x$ in the corresponding Lagrange interpolation polynomials.

The values of the coefficients $\alpha_n(m)$ for different $n$ are given in Table~\ref{tab:alpha}. For example, if $n = 1$, we choose the linear interpolation, which means the rough approximation. Table~\ref{tab:alpha} shows that the coefficients $\alpha_n(m)$ decrease rapidly with increasing $m$.

\begin{table}
  \centering
  \begin{tabular}{|c|c|c|c|}
    \hline
    $n$ & $m=1$ & $m=2$ & $m=3$ \\
    \hline
    $1$ & $0.5000$ & $0$ & $0$ \\
    \hline
    $2$ & $0.5625$ & $-0.0208$ & $0$ \\
    \hline
    $3$ & $0.5859$ & $-0.0326$ & $0.0023$ \\
    \hline
  \end{tabular}
  \caption{The values of coefficient $\alpha_n(m)$ in Eq.~\eqref{eq:firstder}
  for different $n$ and $m$.  
  \label{tab:alpha}}
\end{table}

As one can see in Eq.~\eqref{eq:firstder}, the procedure for finding the first derivative of a discrete function using the weight coefficients $\alpha_n(m)$ is actually a convolution of $S(x)$ and $\alpha_n(m)$, i.e., the cross-correlation function. However, it is known that the cross-correlation function can be defined as the inverse Fourier transform from the product of the spectrum of the first function to the complex conjugate spectrum of the second one~\citep{RabSch07}. Hence it is convenient to analyze the spectral characteristics of both the original discrete function $S(x)$ and the coefficients $\alpha_n(m)$.

\begin{figure}
  \centering
  \includegraphics[scale=0.65]{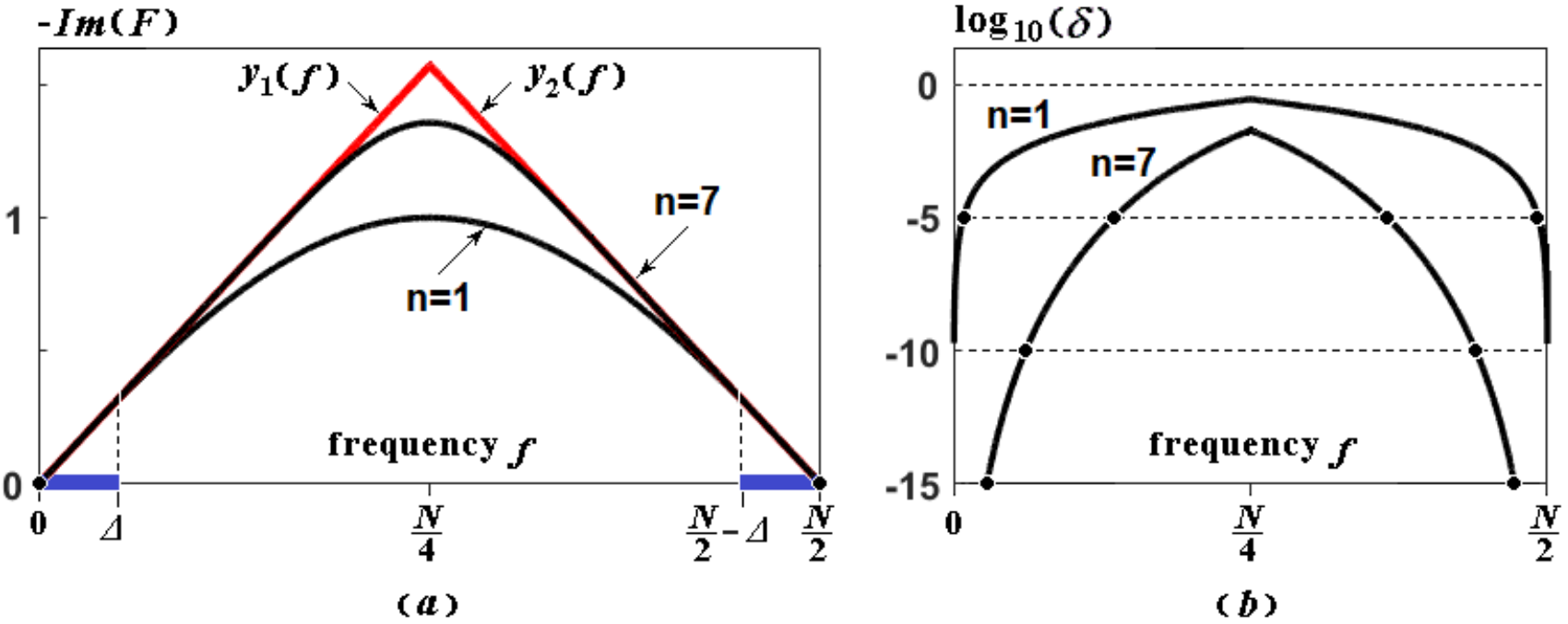}
  \caption{(a) The imaginary part of the spectra of the differential filters $\alpha_n(m)$
  for $n = 1$ and $n = 7$. (b) Errors of the filters.  
  \label{fig:alpha0N/2}}
\end{figure}

It follows from the direct discrete Fourier transform that, for example, for $n = 1$, the spectrum of coefficients $\alpha_1(m)$ has the form,
\begin{equation*}
  F(f) = \frac{1}{2}
  \left[
    \exp
    \left(
      -\mathrm{i}\frac{2\pi}{N}
    \right) - 
    \exp
    \left(
      +\mathrm{i}\frac{2\pi}{N}
    \right)
  \right] =
  - \mathrm{i} \sin
  \left(
    \frac{2\pi}{N} f
  \right).
\end{equation*}
The imaginary parts of the spectral coefficients for $n = 1$ and $n = 7$ are shown in Fig.~\ref{fig:alpha0N/2}(a).

Note that these filters perform the differentiation in two frequency regions. It can be seen that, for $n = 1$, the spectrum is close to the spectra of the limiting differentiating filters $y_1 = 2\pi f/N$ and $y_2 = \pi - 2\pi f/N$, shown by straight red lines, only in narrow bands of spatial frequencies: the zero frequency $[0]$ and the maximal one $[N/2]$. One can see in Fig.~\ref{fig:alpha0N/2}(b) that it is possible to estimate the frequency dependence of the deviation of the spectra of the corresponding coefficients $\alpha_n(m)$ from the spectra of the limiting filters $y_1$ and $y_2$.

To construct one dimensional wave structures, we perform the iteration procedure with respect to the parameter $\tau$,
\begin{equation*}
  S(x,\tau+1) = S(x,\tau) +
  \sum_{m=1}^n
  \left[
    S(x+2m-1) - S(x-2m+1)
  \right]
  \alpha_n(m),
\end{equation*}
where $\tau = 0,1,2,\dotsc$. In order to reduce the complexity of the expression, we will use only the differentiating filter $\alpha_1(m)$,
\begin{equation}\label{eq:S1}
  S(x,\tau+1) = S(x,\tau) +
  \frac{1}{2}
  \left[
    S(x+1) - S(x-1)
  \right].
\end{equation}
Using the representation for $S$ in Eq.~\eqref{eq:cos2sin2}, it can be shown that procedure in Eq.~\eqref{eq:S1} is equivalent to the following approximate differential relations:
\begin{equation*}
  S(x,\tau+1) - S(x,\tau) \approx \frac{\partial p}{\partial \tau},
  \quad
   \frac{1}{2}
  \left[
    S(x+1) - S(x-1)
  \right] \approx
  \frac{\partial q}{\partial x},
\end{equation*}
for even $x$, and
\begin{equation*}
  S(x,\tau+1) - S(x,\tau) \approx \frac{\partial q}{\partial \tau},
  \quad
  \frac{1}{2}
  \left[
    S(x+1) - S(x-1)
  \right] \approx
  \frac{\partial p}{\partial x},
\end{equation*}
for odd $x$. Combining these expressions, we obtain the system of approximate wave differential equations
\begin{equation}\label{eq:sys0N/2}
  \frac{\partial p}{\partial \tau} \approx \frac{\partial q}{\partial x},
  \quad
  \frac{\partial q}{\partial \tau} \approx \frac{\partial p}{\partial x}.
\end{equation}
One can see that the system in Eq.~\eqref{eq:sys0N/2} is analogous to that in Eq.~\eqref{eq:gensys}.

\begin{figure}
  \centering
  \includegraphics[scale=0.7]{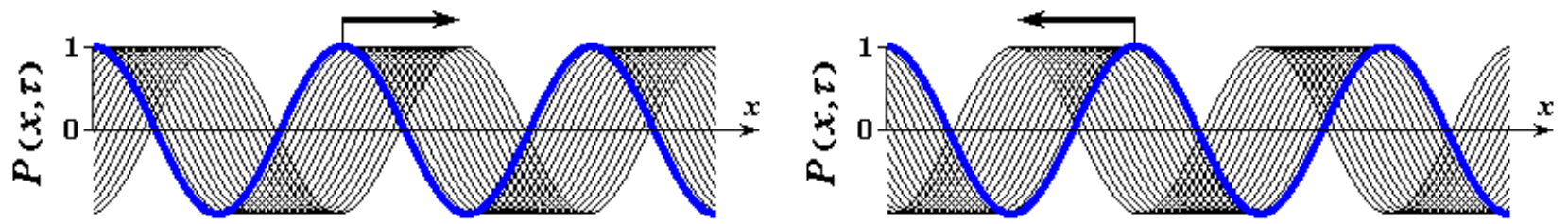}
  \caption{An example of the evolution of a harmonic function corresponding
  to the numerical simulation by the algorithm in Eq.~\eqref{eq:S1}.  
  \label{fig:runningwaves}}
\end{figure}

An example of the numerical simulation using the algorithm in Eq.~\eqref{eq:S1} for harmonic functions is given in Fig.~\ref{fig:runningwaves}.  Here, the initial states of one of the discrete functions are shown by blue lines. If the spectral component of the function $S$ is in the spatial frequency band $[0]$, then it moves to the right. When it is in the frequency band $[N/2]$, the motion is to the left.

As an initial condition of the iterative procedure, one can use a unity, placed in an arbitrary point of the one dimensional grid. At the same time, for all spatial frequencies the absolute value of the initial state spectrum will also be equal to one. We call this initial condition as the shock excitation of the grid. In this case, the two dimensional spectrum of the function $S(x,\tau)$ obtained by the iterative procedure in Eq.~\eqref{eq:S1} can be considered as the amplitude-frequency characteristic (AFC) of the algorithm.

\begin{figure}
  \centering
  \includegraphics[scale=0.7]{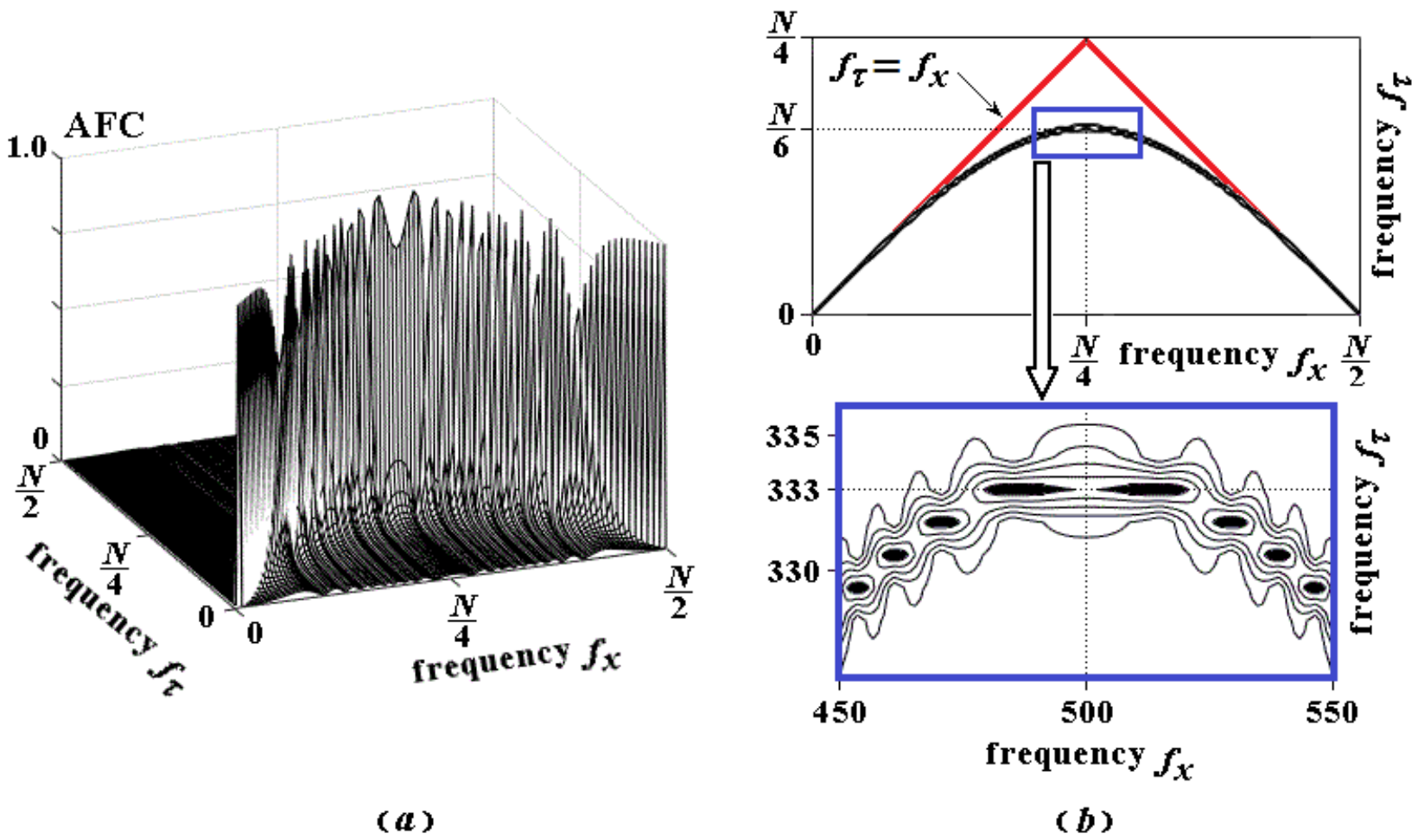}
  \caption{(a) Two-dimensional AFC for the iterative algorithm in Eq.~\eqref{eq:S1} at $n = 1$. 
  (b) The horizontal cross-sections of AFC.
  \label{fig:AFC0N/2}}
\end{figure}

The numerical simulation of the iteration procedure in Eq.~\eqref{eq:S1} was carried out for a one dimensional grid with $N = 2000$ elements, with help of the differentiating filter of the order $n = 1$ with the shock initial condition. The number of iterations is $K =2000$. Using the obtained values   for each number of the interation $\tau = k$, we formed the two dimensional array  $S(x,k)$ having the size $(N \times K) = (2000 \times 2000)$. In Fig.~\ref{fig:AFC0N/2}(a), we show the two dimensional spectrum of   in the one dimensional grid with $N = 2000$ elements and for $n = 1$. In Fig.~\ref{fig:AFC0N/2}(b), the horizontal cross-sections are given.

From the given data, it can be seen that, in the narrow band of spatial frequencies near zero, the following condition is met: $f_\tau \approx f_x$. Thus the group velocity reads $\mathrm{d}f_\tau / \mathrm{d} f_x \approx 1$. For the maximal frequency $N/2$ one gets that $f_\tau \approx N/2 - f_x$ and $\mathrm{d}f_\tau / \mathrm{d} f_x \approx -1$. Thus, it can be noted that for these bands of spatial frequencies the algorithm in Eq.~\eqref{eq:S1} satisfactorily generates one dimensional linear wave structures. Note that, for other frequency bands, there is a more complex behavior. For example, at the bottom of Fig.~\ref{fig:AFC0N/2}(b) an enlarged fragment of the cross-section for the central spatial frequency $N/4$ is shown. However, the explanation of this feature is beyond the scope of this work.

\section{One dimensional grid: Frequency band $[N/4]$}\label{sec:1DN/4}

\begin{figure}
  \centering
  \includegraphics[scale=0.6]{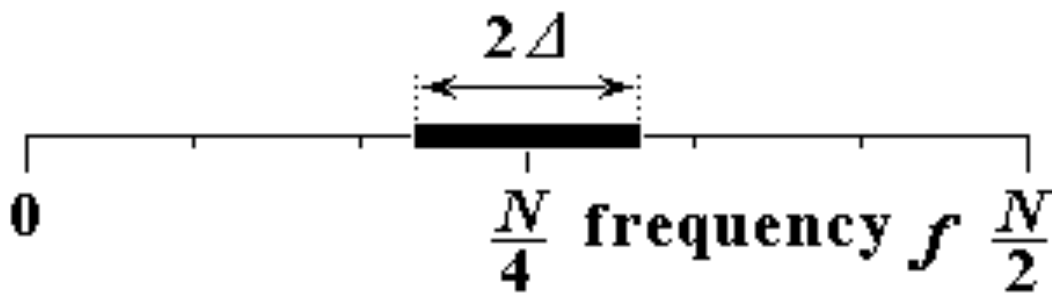}
  \caption{Frequency band of the discrete function $S(x)$.
  \label{fig:freqN/4}}
\end{figure}

Let us now suppose that the spectrum of the one dimensional discrete function $S(x)$ is located near the central spatial frequency $N/4$, as shown by the bold segment in Fig.~\ref{fig:freqN/4}. We denote this band as $[N/4]$. Using the inverse discrete Fourier transform, one can show that, in this case, $S(x)$ can be represented in the form,
\begin{equation}\label{eq:cossin}
  S(x) = p(x) \cos \frac{\pi x}{2} + q(x) \sin \frac{\pi x}{2},
\end{equation}
One can see that for even $x$, $S(x) = \pm p$, whereas for odd $x$: $S(x) = \pm q$. Let us call $p$ and $q$ virtual sign alternating functions.

\begin{figure}
  \centering
  \includegraphics[scale=0.45]{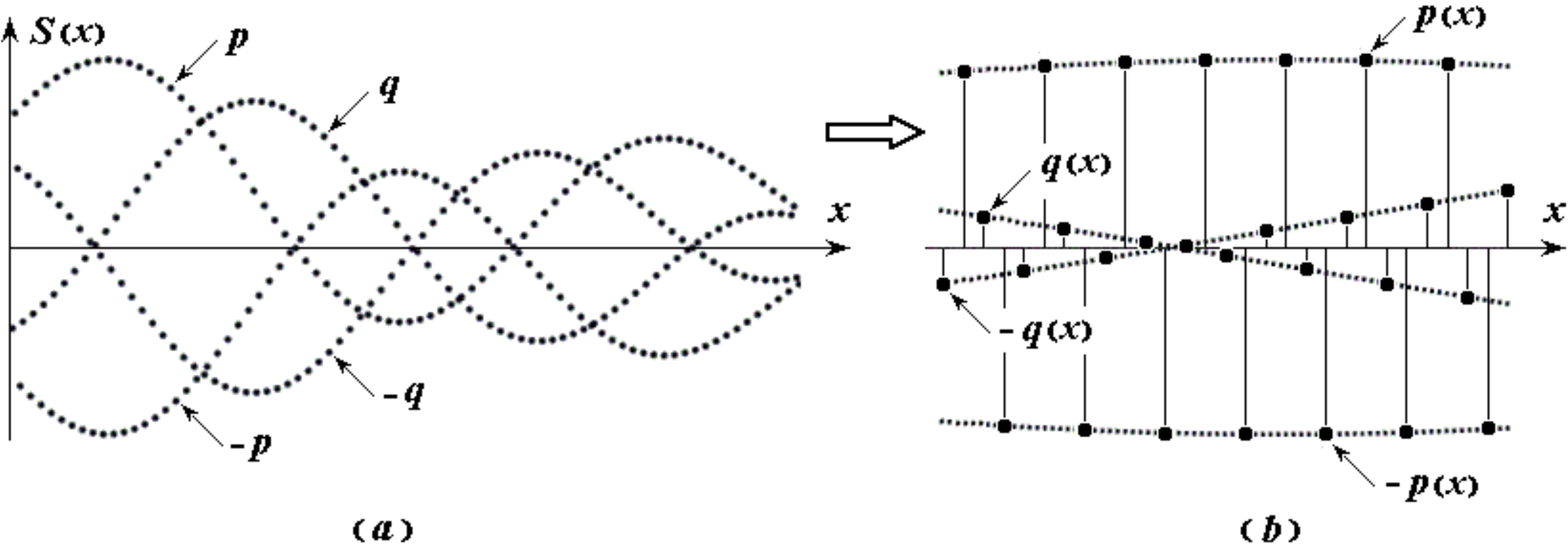}
  \caption{(a) The sample of the discrete function for the frequency band $[N/4]$.
  (b) The fragment of such a function.
  \label{fig:sampleN/4}}
\end{figure}

In Fig.~\ref{fig:sampleN/4}(a), we show an example of a discrete function in this frequency band, assuming the condition that $\Delta \ll N$. A fragment of such a function is depicted in Fig.~\ref{fig:sampleN/4}(b). Fig.~\ref{fig:sampleN/4} shows that the original discontinuous function $S(x)$ can be represented as the sum of two relatively ``smooth'' sign alternating functions $\pm p$ and $\pm q$, where $p$ is defined in even $x$, and $q$ in odd $x$.

The expression for the weight coefficients $\alpha_n(m)$ in the considered frequency band was obtained by~\cite{Dvo07,Dvo08},
\begin{equation}\label{eq:alphaN/4}
  \alpha_n(m) =
  \left|
    2(2m-1)
    \prod_{\substack{k=1 \\ k\neq m}}^{n}
    \left[
      1 - \frac{(2m-1)^2}{(2m-1)^2}
    \right]
  \right|^{-1},
\end{equation}
where $m = 1,2, \dots n$. The values of the coefficients $\alpha_n(m)$ for several $n$ are given in Table~\ref{tab:alphaN/4}.

\begin{table}
  \centering
  \begin{tabular}{|c|c|c|c|}
    \hline
    $n$ & $m=1$ & $m=2$ & $m=3$ \\
    \hline
    $1$ & $0.5000$ & $0$ & $0$ \\
    \hline
    $2$ & $0.5625$ & $0.0208$ & $0$ \\
    \hline
    $3$ & $0.5859$ & $0.0326$ & $0.0023$ \\
    \hline
  \end{tabular}
  \caption{The values of coefficient $\alpha_n(m)$ in Eq.~\eqref{eq:alphaN/4}
  for different $n$ and $m$.  
  \label{tab:alphaN/4}}
\end{table}

For example, the situation when $n = 1$ corresponds to the linear interpolation and a rough approximation. Table~\ref{tab:alphaN/4} shows that the coefficients $\alpha_n(m)$ decrease rapidly with increasing $m$. Using the discrete Fourier transform, one gets that the spectrum of coefficients $\alpha_n(m)$, e.g., for $n = 1$, has the form,
\begin{equation*}
  F(f) = \frac{1}{2}
  \left[
    \exp
    \left(
      -\mathrm{i}\frac{2\pi}{N}
    \right) + 
    \exp
    \left(
      +\mathrm{i}\frac{2\pi}{N}
    \right)
  \right] =
  \cos
  \left(
    \frac{2\pi}{N} f
  \right).
\end{equation*}

\begin{figure}
  \centering
  \includegraphics[scale=0.6]{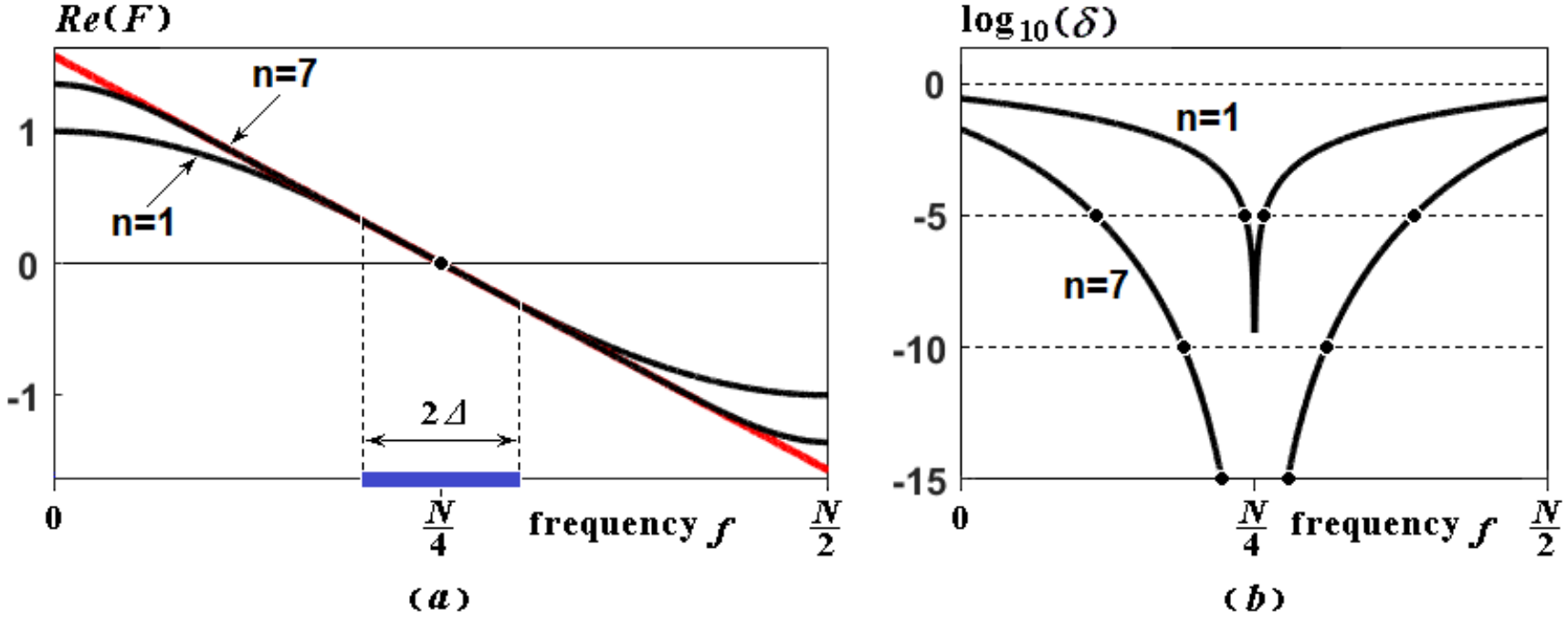}
  \caption{(a) The real part of the spectra of the differential filters $\alpha_n(m)$
  for $n = 1$ and $n = 7$.
  (b) The errors of the filters.
  \label{fig:specN/4}}
\end{figure}

The real parts of the spectra of the coefficients $\alpha_n(m)$ for the values $n = 1$ and $n = 7$ are shown in Fig.~\ref{fig:specN/4}(a). Thus, for $n = 1$, the part of the spectrum of $\alpha_n(m)$, which is in a narrow band of spatial frequencies near the central frequency $N/4$, is close to the spectrum of the limiting differentiating filter, which is also shown in Fig.~\ref{fig:specN/4}(a) with the red line. The deviation of the spectra of the corresponding coefficients $\alpha_n(m)$ from the spectrum of the limiting filter can be estimated using Fig.~\ref{fig:specN/4}(b).

To construct one dimensional wave structures in the considered frequency band, it is necessary to implement the following iteration procedure with respect to the parameter $\tau$:
\begin{equation}\label{eq:S1N/4}
  S(x,\tau+1) = S(x,\tau) +
  \frac{1}{2}
  \left[
    S(x+1) + S(x-1)
  \right] X(x),
\end{equation}
where $\tau = 0,1,2,\dots$ and $X(x) = (-1)^x$ is the additional multiplier.

Using the expression for $S$ in Eq.~\eqref{eq:cossin}, it can be shown that the procedure in Eq.~\eqref{eq:S1N/4} is equivalent to the following approximate differential relations:
\begin{equation*}
  S(x,\tau+1) - S(x,\tau) \approx \frac{\partial p}{\partial \tau},
  \quad
  \frac{1}{2}
  \left[
    S(x+1) + S(x-1)
  \right]X(x) \approx
  \frac{\partial q}{\partial x},
\end{equation*}
for even $x$, and
\begin{equation*}
  S(x,\tau+1) - S(x,\tau) \approx \frac{\partial q}{\partial \tau},
  \quad
  \frac{1}{2}
  \left[
    S(x+1) + S(x-1)
  \right]X(x) \approx
  \frac{\partial p}{\partial x},
\end{equation*}
for odd $x$. Combining these expressions, we obtain the system of approximate wave differential equations
\begin{equation}\label{eq:sys0N/4}
  \frac{\partial p}{\partial \tau} \approx \frac{\partial q}{\partial x},
  \quad
  \frac{\partial q}{\partial \tau} \approx \frac{\partial p}{\partial x},
\end{equation}
which again coincide with Eq.~\eqref{eq:gensys}.

The iterative procedure for the frequency band $[N/4]$ in Eq.~\eqref{eq:S1N/4} differs significantly from that in Eq.~\eqref{eq:S1}. Here, to obtain the approximate differential Eq.~\eqref{eq:sys0N/4}, the equidistant values of $S$ are summed and an additional multiplier $X(x)$ is introduced for the implementation of the system leading to the wave propagation.

\begin{figure}
  \centering
  \includegraphics[scale=0.7]{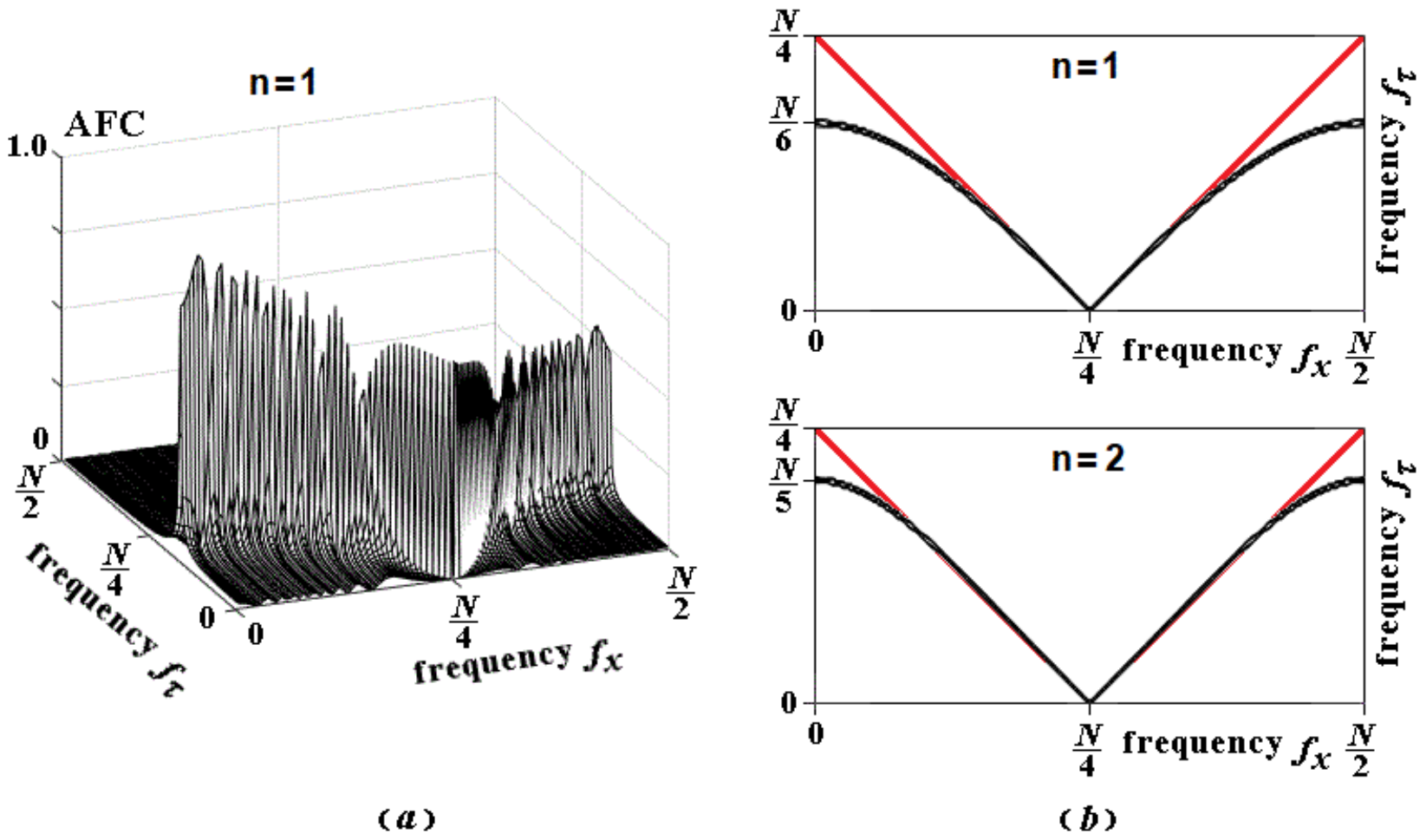}
  \caption{(a) Two dimensional AFC, for $n = 1$. (b) The horizontal cross sections of AFC
  for $n = 1$ and $n = 2$.
  \label{fig:AFCN/4}}
\end{figure}

In Fig.~\ref{fig:AFCN/4}(a), we show the two dimensional spectrum (AFC), obtained using the iterative procedure in Eq.~\eqref{eq:S1N/4}, for $n = 1$ under the shock initial condition. The horizontal cross-sections of the AFC for $n = 1$ and $n = 2$ are depicted in Fig.~\ref{fig:AFCN/4}(b).

One can see in Fig.~\ref{fig:AFCN/4} that, in the narrow band of spatial frequencies below the central frequency $N/4$, the following condition is satisfied: $f_\tau \approx N/4 - f_x$, whereas for frequencies above $N/4$, one has $f_\tau \approx f_x - N/4$. For these bands, the group velocities are approximately equal to $+1$ and $-1$, respectively. It means that in the considered narrow frequency band, the algorithm in Eq.~\eqref{eq:S1N/4} leads to the formation of one dimensional linear wave structures. In Fig.~\ref{fig:AFCN/4}(a), one can observe some asymmetry of the two dimensional spectrum at large deviations of $f_x$ from the spatial frequency $N/4$. The analysis of this feature of the spectrum will be given in one of our forthcoming works.

\section{Two dimensional grid: Frequency zones $[0] + [N/2]$}\label{sec:2D0N/2}

\begin{figure}
  \centering
  \includegraphics[scale=0.6]{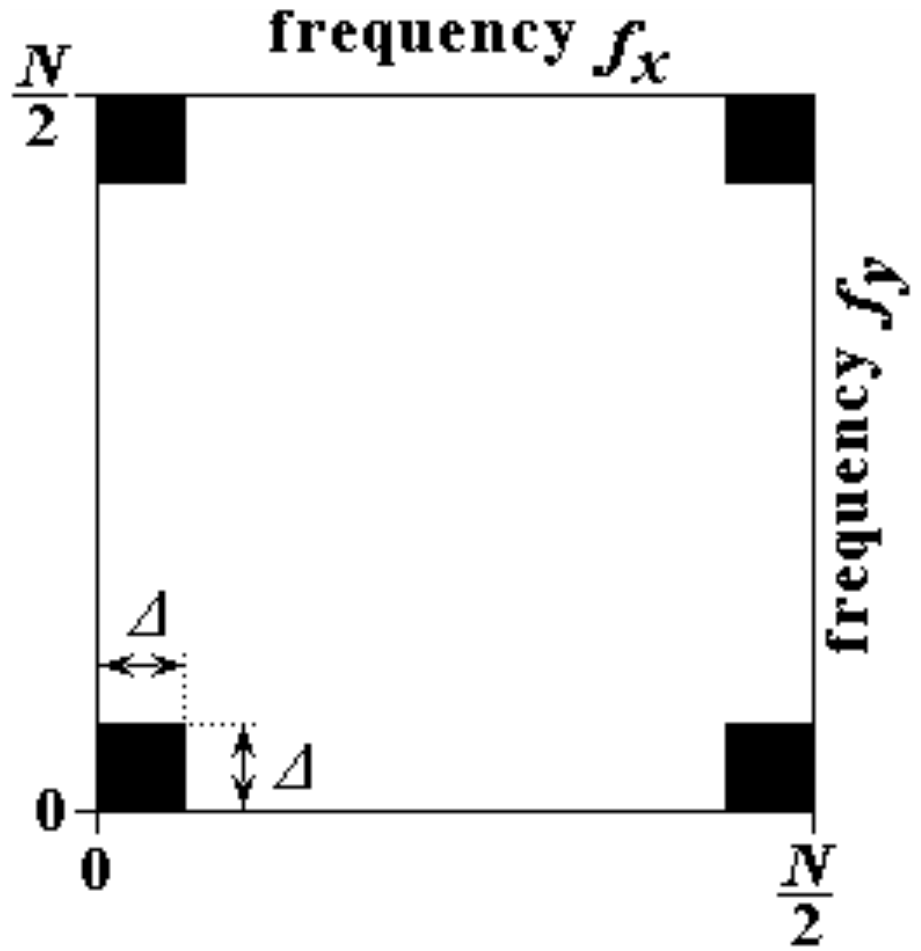}
  \caption{Frequency zones $[0] + [N/2]$.
  \label{fig:2Dfreq0N/2}}
\end{figure}

Now we turn to the studies of two dimensional grids. If the spectrum of the two dimensional discrete function $S(x,y)$ lies in the frequency bands shown in Fig.~\ref{fig:2Dfreq0N/2}, then, using the inverse two dimensional Fourier transform, we can show that this function can be represented in the form,
\begin{align}\label{eq:2Dcos2sin2}
  S(x,y) = & p \cos^2 \frac{\pi x}{2} \cos^2 \frac{\pi y}{2} +
  r \sin^2 \frac{\pi x}{2} \cos^2 \frac{\pi y}{2}
  \notag
  \\
  & +
  q \sin^2 \frac{\pi x}{2} \sin^2 \frac{\pi y}{2} +
  s \cos^2 \frac{\pi x}{2} \sin^2 \frac{\pi y}{2},
\end{align}
where $x,y = 0,1,\dots, N-1$.

\begin{figure}
  \centering
  \includegraphics[scale=0.8]{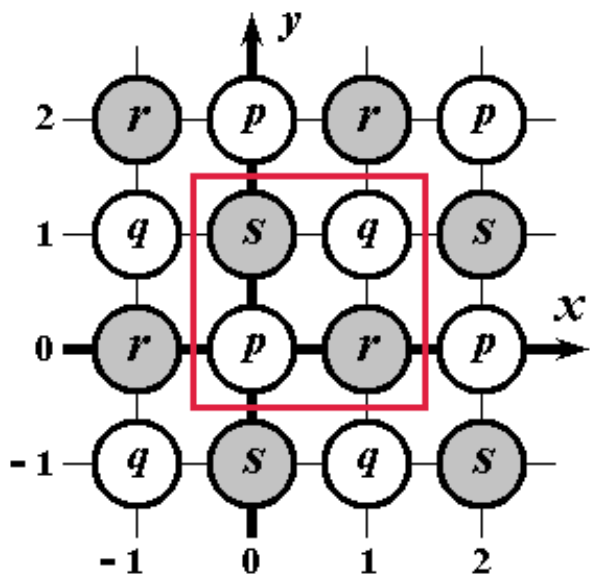}
  \caption{The location of 4 virtual functions for a two dimensional grid.
  The size of the fragment is $4\times 4$.
  \label{fig:2Dvirtfun0N/2}}
\end{figure}

Using Eq.~\eqref{eq:2Dcos2sin2}, we can verify that, depending on whether $x$ and $y$ are even or odd, we always have only one of the four functions $p$, $q$, $r$, and $s$, as shown in Fig.~\eqref{fig:2Dvirtfun0N/2}. Thus, for $\Delta \ll N$, the original discontinuous function $S(x,y)$ can be represented as the sum of four relatively smooth, linearly independent virtual functions.

Then, we implement the iteration procedure with respect to the parameter $\tau$ as
\begin{align}\label{eq:2DS1}
  S(x,y,\tau+1) = & S(x,y,\tau) +
  \frac{1}{2}
  \left[
    S(x+1,y) - S(x-1,y)
  \right]X(x,y) 
  \notag
  \\
  & +
  \frac{1}{2}
  \left[
    S(x,y+1) - S(x,y-1)
  \right]Y(x,y),
\end{align}
where $\tau = 0,1,2,\dotsc$. The multipliers $X(x,y)=1$ and $X(x,y)=(-1)^x$ in Eq.~\eqref{eq:2DS1} are necessary for the implementation of the wave equation (see below).

Using Eq.~\eqref{eq:2Dcos2sin2}, we can show that procedure in Eq.~\eqref{eq:2DS1} is equivalent to approximate differential relations,
\begin{equation}\label{eq:2Dsys0N/2}
  \frac{\partial p}{\partial \tau} \approx \frac{\partial r}{\partial x} - \frac{\partial s}{\partial y},
  \quad
  \frac{\partial q}{\partial \tau} \approx \frac{\partial s}{\partial x} + \frac{\partial r}{\partial y},
  \quad
  \frac{\partial r}{\partial \tau} \approx \frac{\partial p}{\partial x} + \frac{\partial q}{\partial y},
  \quad
  \frac{\partial s}{\partial \tau} \approx \frac{\partial q}{\partial x} - \frac{\partial p}{\partial y}.
\end{equation}
Eq.~\eqref{eq:2Dsys0N/2} is identical to the system of two complex differential equations,
\begin{equation}\label{eq:2Dcomp0N/2}
  \frac{\partial Z_1}{\partial \tau} = D_2 Z_2,
  \quad
  \frac{\partial Z_2}{\partial \tau} = D_2^* Z_1,
\end{equation}
where $Z_1 = p + \mathrm{i}q$, $Z_2 = r + \mathrm{i}s$, and $D_2 = \tfrac{\partial}{\partial x} + \mathrm{i} \tfrac{\partial}{\partial y}$. The solution of Eq.~\eqref{eq:2Dcomp0N/2} has the form,
\begin{align}\label{eq:2Dcomp0N/2}
  Z_1 = & Z_1^{(0)} \exp[2\pi\mathrm{i}(f_\tau \tau + f_x x + f_y y)],
  \notag
  \\
  Z_2 = & Z_2^{(0)} \exp[2\pi\mathrm{i}(f_\tau \tau + f_x x + f_y y)],
\end{align}
where $Z_{1,2}^{(0)}$ are some constant complex numbers and $f_\tau^2 = f_x^2 + f_y^2$.

The numerical simulation of the iterative procedure in Eq.~\eqref{eq:2DS1} is carried out for the two dimensional grid with the size $(N \times N) = (400 \times 400)$, for the differentiating filter of the order $n = 2$ at the shock initial condition. The number of iterations is $K=400$. Using the obtained values $S(x,y,k)$ for each number of the iteration $\tau = k$, we form the three dimensional array $S(x,y,\tau)$ of the size $(N \times N \times K) = (400 \times 400 \times 400)$. Then we compute the three dimensional spectrum of the function $S(x,y,\tau)$ and find the coordinates $(f_x, f_y, f_\tau)$ for all values of this spectrum which are above a certain threshold value.

\begin{figure}
  \centering
  \includegraphics[scale=0.8]{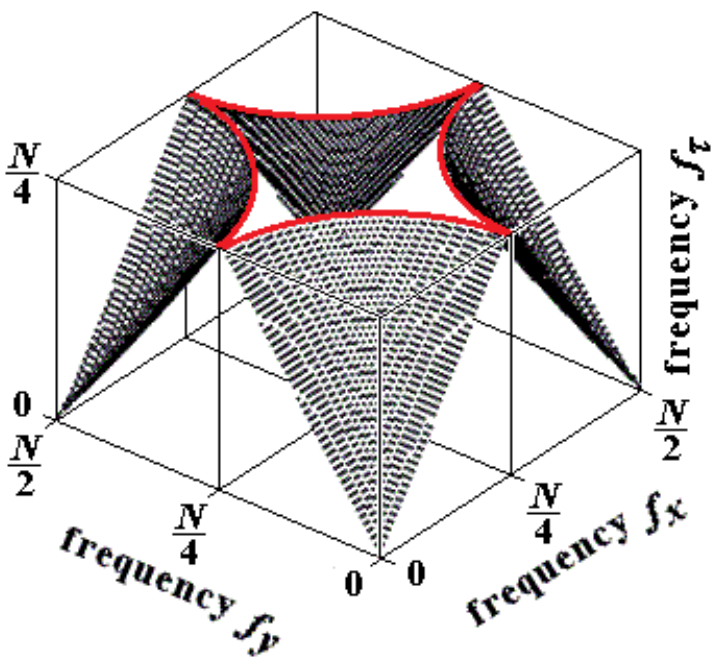}
  \caption{The location of the maximum values of the three dimensional spectrum of the function
  $S(x, y, \tau)$.
  \label{fig:2Dmax0N/2}}
\end{figure}

An example of the numerical simulation by the algorithm in Eq.~\eqref{eq:2DS1} is shown in Fig.~\ref{fig:2Dmax0N/2}. For the obtained function $S(x,y,\tau)$, the maximal values of its three dimensional spectrum are located on the surfaces of four cones. For $\Delta \ll N$, the absolute value of the group velocities of waves is approximately equal to one:
\begin{equation*}
  \sqrt{\left( \frac{\partial f_\tau}{\partial f_x} \right)^2 +
  \left( \frac{\partial f_\tau}{\partial f_y} \right)^2} \approx 1.
\end{equation*}
This fact confirms that the solutions obtained are wave structures.

\section{Two dimensional grid: Frequency zone $[N/4]$}\label{sec:2DN/4}

\begin{figure}
  \centering
  \includegraphics[scale=0.6]{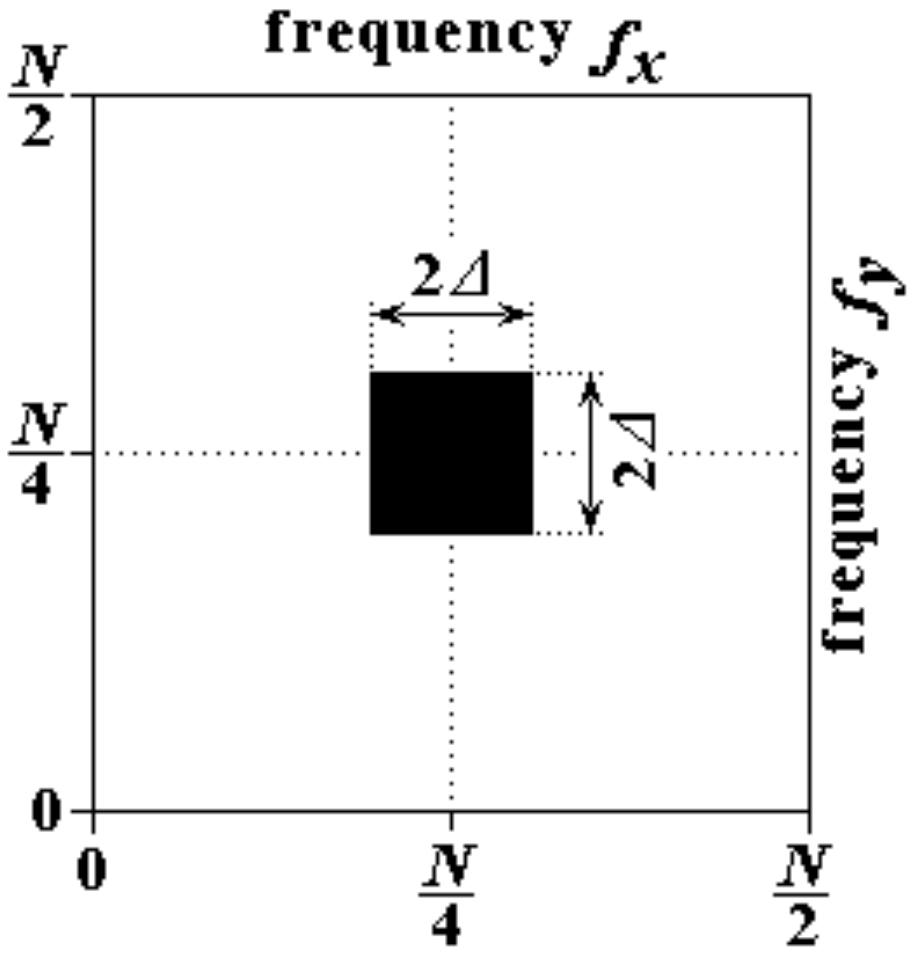}
  \caption{Frequency zone $[N/4]$.
  \label{fig:2DfreqN/4}}
\end{figure}

If the spectrum of the two dimensional discrete function $S(x, y)$ is located in the frequency band shown in Fig.~\ref{fig:2DfreqN/4}, then this function can be represented in the form,
\begin{align}\label{eq:2Dcossin}
  S(x,y) = & p \cos \frac{\pi x}{2} \cos \frac{\pi y}{2} +
  r \sin \frac{\pi x}{2} \cos \frac{\pi y}{2}
  \notag
  \\
  & +
  q \sin \frac{\pi x}{2} \sin \frac{\pi y}{2} +
  s \cos \frac{\pi x}{2} \sin \frac{\pi y}{2},
\end{align}
where $x,y = 0,1,\dots, N-1$.

\begin{figure}
  \centering
  \includegraphics[scale=0.8]{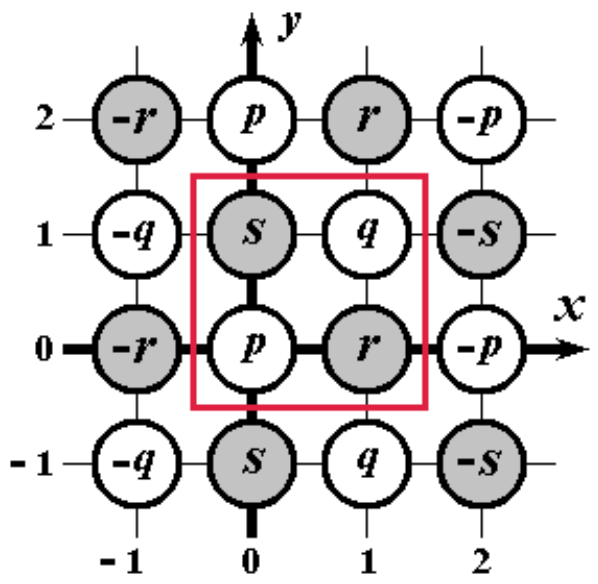}
  \caption{The location of 4 virtual functions for a two dimensional grid.
  The size of the fragment is $4\times 4$.
  \label{fig:2DvirtfunN/4}}
\end{figure}

Using Eq.~\eqref{eq:2Dcossin}, we can verify that, depending whether $x$ and $y$ are even or odd, we will always have only one of the four functions $\pm p$, $\pm q$, $\pm r$, and $\pm s$, as shown in Fig.~\ref{fig:2DvirtfunN/4}. For $\Delta \ll N$, the initial discontinuous function $S(x, y)$ can be represented as the sum of four relatively smooth, virtual sign alternating functions.

The iteration procedure for this case has the form,
\begin{align}\label{eq:2DSN/4}
  S(x,y,\tau+1) = & S(x,y,\tau) +
  \frac{1}{2}
  \left[
    S(x+1,y) + S(x-1,y)
  \right]X(x,y) 
  \notag
  \\
  & +
  \frac{1}{2}
  \left[
    S(x,y+1) + S(x,y-1)
  \right]Y(x,y),
\end{align}
where $\tau = 0,1,2,\dotsc$. The additional multipliers are $X(x,y)=(-1)^x$ and $X(x,y)=(-1)^{x+y}$. Using Eq.~\eqref{eq:2Dcossin}, it can be shown that the procedure in Eq.~\eqref{eq:2DSN/4} is also equivalent to approximate differential relations in Eq.~\eqref{eq:2Dsys0N/2}.

\begin{figure}
  \centering
  \includegraphics[scale=0.8]{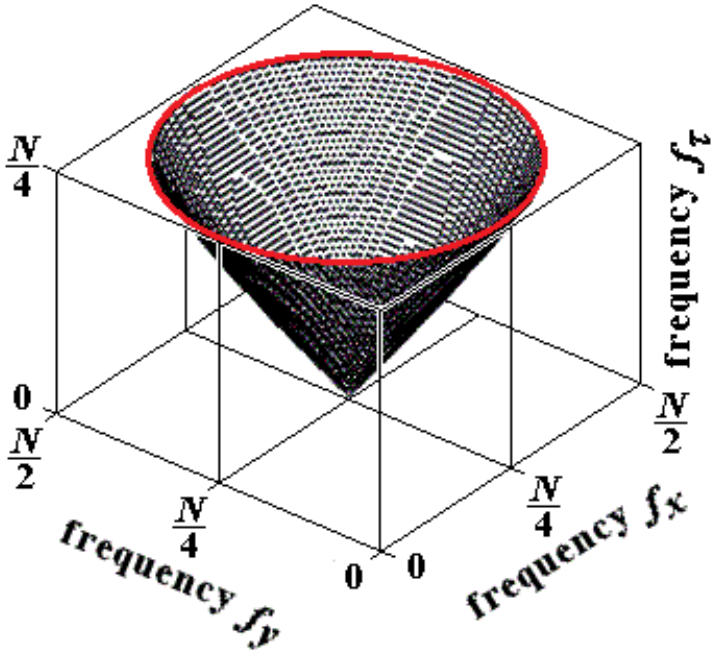}
  \caption{The location of the maximum values of the three dimensional spectrum of the function
  $S(x, y, \tau)$.
  \label{fig:2DmaxN/4}}
\end{figure}

An example of the numerical simulation by the algorithm in Eq.~\eqref{eq:2DSN/4} is given in Fig.~\ref{fig:2DmaxN/4}. It is seen that the obtained function $S(x, y, \tau)$ has the maximum values of its three dimensional spectrum located on the cone surface. If $\Delta \ll N$, the absolute values of the group velocities of waves are approximately equal to one:
\begin{equation*}
  \sqrt{\left( \frac{\partial f_\tau}{\partial f_x} \right)^2 +
  \left( \frac{\partial f_\tau}{\partial f_y} \right)^2} \approx 1.
\end{equation*}
It also confirms that the obtained solutions are wave structures.

\section{Three dimensional grid}\label{sec:3D}

\begin{figure}
  \centering
  \includegraphics[scale=0.8]{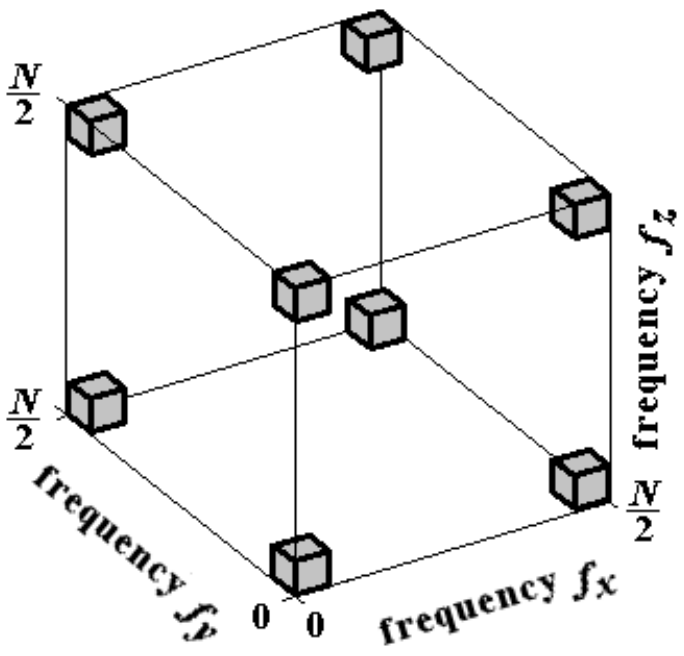}
  \caption{Frequency zones $[0] + [N/2]$.
  \label{fig:3Dfreq0N/2}}
\end{figure}

If the spectrum of the three dimensional discrete function $S(x, y, z)$ is located in the frequency zones shown in Fig.~\ref{fig:3Dfreq0N/2}, then this function can be represented in the form,
\begin{align}\label{eq:3Dcos2sin2}
  S(x,y,z) = &
  p_0 \cos^2 \frac{\pi x}{2} \cos^2 \frac{\pi y}{2} \cos^2 \frac{\pi z}{2} +
  q_0 \sin^2 \frac{\pi x}{2} \sin^2 \frac{\pi y}{2} \sin^2 \frac{\pi z}{2}
  \notag
  \\
  & +
  p_1 \cos^2 \frac{\pi x}{2} \sin^2 \frac{\pi y}{2} \sin^2 \frac{\pi z}{2} +
  q_1 \sin^2 \frac{\pi x}{2} \cos^2 \frac{\pi y}{2} \cos^2 \frac{\pi z}{2}
  \notag
  \\
  & +
  p_2 \sin^2 \frac{\pi x}{2} \cos^2 \frac{\pi y}{2} \sin^2 \frac{\pi z}{2} +
  q_2 \cos^2 \frac{\pi x}{2} \sin^2 \frac{\pi y}{2} \cos^2 \frac{\pi z}{2}
  \notag
  \\
  & +
  p_3 \sin^2 \frac{\pi x}{2} \sin^2 \frac{\pi y}{2} \cos^2 \frac{\pi z}{2} +
  q_3 \cos^2 \frac{\pi x}{2} \cos^2 \frac{\pi y}{2} \sin^2 \frac{\pi z}{2}.  
\end{align}
Using Eq.~\eqref{eq:3Dcos2sin2}, we obtain that, depending whether $x$, $y$ and $z$ are even or odd, we always have only one of the eight functions $p_0, \dots p_3$ and $q_0, \dots q_3$. For $\Delta \ll N$ the original discontinuous function $S(x, y, z)$ can be represented as the sum of eight relatively smooth, linearly independent virtual functions.

The iterative procedure in this case has the form,
\begin{align}\label{eq:3DS1}
  S(x,y,z,\tau+1) = & S(x,y,z,\tau) +
  \frac{1}{2}
  \left[
    S(x+1,y,z) - S(x-1,y,z)
  \right]X(x,y,z) 
  \notag
  \\
  & +
  \frac{1}{2}
  \left[
    S(x,y+1,z) - S(x,y-1,z)
  \right]Y(x,y,z)
  \notag
  \\
  & +
  \frac{1}{2}
  \left[
    S(x,y,z+1) - S(x,y,z-1)
  \right]Z(x,y,z),
\end{align}
where $\tau = 0, 1, 2, \dotsc$. The additional multipliers are $X(x,y,z)=(-1)^y$,  $Y(x,y,z)=(-1)^z$, and $Z(x,y,z)=(-1)^x$.

Using the representation for $S(x, y, z)$ in Eq.~\eqref{eq:3Dcos2sin2}, we can show that procedure in Eq.~\eqref{eq:3DS1} is equivalent to the system of approximate differential relations,
\begin{align}\label{eq:3Dsys0N/2}
  \frac{\partial p_0}{\partial \tau} \approx &
  - \frac{\partial q_1}{\partial x} - \frac{\partial q_2}{\partial y} - \frac{\partial q_3}{\partial z},
  \quad
  \frac{\partial p_1}{\partial \tau} \approx
  \frac{\partial q_0}{\partial x} + \frac{\partial q_3}{\partial y} - \frac{\partial q_2}{\partial z},
  \notag
  \\
  \frac{\partial p_2}{\partial \tau} \approx &
  - \frac{\partial q_3}{\partial x} + \frac{\partial q_0}{\partial y} + \frac{\partial q_1}{\partial z},
  \quad
  \frac{\partial p_3}{\partial \tau} \approx 
  \frac{\partial q_2}{\partial x} - \frac{\partial q_1}{\partial y} + \frac{\partial q_0}{\partial z},
  \notag
  \\
  \frac{\partial q_0}{\partial \tau} \approx &
  \frac{\partial p_1}{\partial x} + \frac{\partial p_2}{\partial y} + \frac{\partial p_3}{\partial z},
  \quad
  \frac{\partial q_1}{\partial \tau} \approx
  - \frac{\partial p_0}{\partial x} - \frac{\partial p_3}{\partial y} + \frac{\partial p_2}{\partial z},
  \notag
  \\
  \frac{\partial q_2}{\partial \tau} \approx &
  \frac{\partial p_3}{\partial x} - \frac{\partial p_0}{\partial y} - \frac{\partial p_1}{\partial z},
  \quad
  \frac{\partial q_3}{\partial \tau} \approx 
  - \frac{\partial p_2}{\partial x} + \frac{\partial p_1}{\partial y} - \frac{\partial p_0}{\partial z}.
\end{align}
The system in Eq.~\eqref{eq:3Dsys0N/2} can be rewritten in the form of two differential quaternion equations:
\begin{equation}\label{eq:3Dcomp0N/2}
  \frac{\partial Z_1}{\partial \tau} = D_3 Z_2,
  \quad
  \frac{\partial Z_2}{\partial \tau} = D_3^* Z_1,
\end{equation}
where $Z_1 = p_0 + \mathrm{i}p_1 + \mathrm{j}p_2 + \mathrm{k}p_3$, $Z_2 = q_0 + \mathrm{i}q_1 + \mathrm{j}q_2 + \mathrm{k}q_3$, and $D_3 = \mathrm{i}\tfrac{\partial}{\partial x} + \mathrm{j} \tfrac{\partial}{\partial y} + \mathrm{k} \tfrac{\partial}{\partial z}$. Here we use the quaternion algebra reviewed by~\cite{Gir84}.

The solution of Eq.~\eqref{eq:3Dcomp0N/2} has the form,
\begin{align}\label{eq:2Dcomp0N/2}
  Z_1 = & Z_1^{(0)} \exp[2\pi\mathrm{i}(f_\tau \tau + f_x x + f_y y + f_z z)],
  \notag
  \\
  Z_2 = & Z_2^{(0)} \exp[2\pi\mathrm{i}(f_\tau \tau + f_x x + f_y y + f_z z)],
\end{align}
where $Z_{1,2}^{(0)}$ are some constant quaternions and $f_\tau^2 = f_x^2 + f_y^2 + f_z^2$.

\begin{figure}
  \centering
  \includegraphics[scale=0.8]{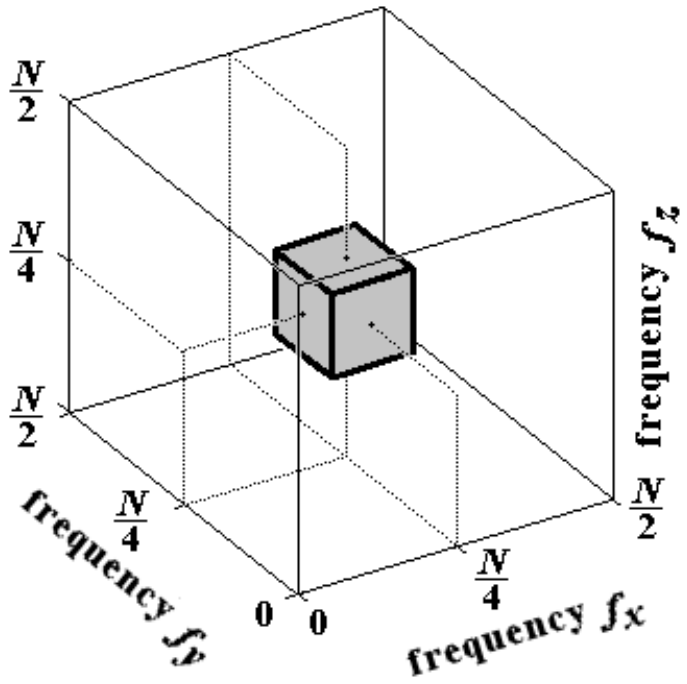}
  \caption{Frequency zone $[N/4]$.
  \label{fig:3DfreqN/4}}
\end{figure}

For a frequency zone $[N/4]$, shown in Fig.~\ref{fig:3DfreqN/4}, the iterative procedure reads
\begin{align}\label{eq:3DSN/4}
  S(x,y,z,\tau+1) = & S(x,y,z,\tau) +
  \frac{1}{2}
  \left[
    S(x+1,y,z) + S(x-1,y,z)
  \right]X(x,y,z) 
  \notag
  \\
  & +
  \frac{1}{2}
  \left[
    S(x,y+1,z) + S(x,y-1,z)
  \right]Y(x,y,z)
  \notag
  \\
  & +
  \frac{1}{2}
  \left[
    S(x,y,z+1) + S(x,y,z-1)
  \right]Z(x,y,z),
\end{align}
where $\tau = 0, 1, 2, \dotsc$. The additional multipliers are $X(x,y,z)=(-1)^{x+y}$,  $Y(x,y,z)=(-1)^{y+z}$, and $Z(x,y,z)=(-1)^{x+z}$.

The limited access to powerful computing facilities does not allow us to perform a detailed analysis of the solutions obtained by iterative algorithms in Eqs.~\eqref{eq:3DS1} and~\eqref{eq:3DSN/4}. The numerical simulation for three dimensional grids of the small size confirms the efficiency of the proposed algorithms. However it does not provide convincing visual results.

\section{Four dimensional grid}\label{sec:4D}

The principle for the construction of linear wave structures in four dimensional grids is identical to the above methods for grids of smaller dimensions. For four dimensional grids the number of discrete virtual functions will be 16 and it is also possible to construct the corresponding iterative procedures for them. Note that it is impossible to construct wave structures in the considered frequency zones in grids with dimensions higher than four.

\section{Conclusion}\label{sec:CONCL}

In the present paper, we have studied the possibility for the construction of linear wave structures in discrete grids with dimensions up to four. We have obtained that, in the numerical construction of linear wave structures, there is no need to assign multiple values to each point of the grid. The consideration of the discrete multidimensional function in the corresponding frequency zones allows one to allocate the required number of virtual functions. To implement the propagation of a wave one has to introduce the additional multipliers. We have also obtained that, using the differentiating filters, one can improve the quality and accuracy of the construction of wave structures. The obtained results can be used in the lattice simulations in quantum field theory.

Some issues have not been considered in our work. In particular, we have hot analyzed methods of filtering in the discrete multidimensional grids. The possibility to use of other frequency zones has not been studied. We have not discussed alternative methods of the iterative procedure implementation, including introduction of nonlinearities. We plan to address these issues in our forthcoming works.

\section*{Acknowledgment}
One of the authors (MD) thanks RFBR (Grant No.~18-02-00149a) for partial support.

\end{document}